\documentclass[11pt,reqno]{amsart}

\usepackage[margin=1.2in]{geometry}
\usepackage{amsmath,amssymb,amsthm}
\usepackage{tikz}
\usepackage{hyperref}

\tikzset{
  vtx/.style={circle,draw,fill=white,inner sep=0pt,minimum size=3.4pt},
  cvtx/.style={circle,draw,fill=black,inner sep=0pt,minimum size=4pt},
  rvtx/.style={circle,draw,fill=gray!55,inner sep=0pt,minimum size=4pt},
  blob/.style={circle,draw,fill=gray!12,inner sep=0pt,minimum size=15pt},
}
\hypersetup{colorlinks=true,linkcolor=blue,citecolor=blue,urlcolor=blue}

\newtheorem{theorem}{Theorem}
\newtheorem{conjecture}[theorem]{Conjecture}

\DeclareMathOperator{\Zf}{Z}

\begin{document}

\title[A counterexample to the $\alpha$--$Z$ conjecture]
{A counterexample to the zero forcing versus independence conjecture
for cubic and subcubic graphs}

\author{Mikko Fischer}
\address{Independent researcher, Raasepori, Finland}
\email{mikko.fischer@gmail.com}
\thanks{The counterexamples were found with the assistance of Claude Opus~5
(Anthropic), directed by the author. The author has independently executed the
verification script included in the ancillary files and checked its output
against the statements made here.}

\date{\today}

\subjclass[2020]{05C69, 05C50}
\keywords{zero forcing number, independence number, cubic graph, TxGraffiti}

\begin{abstract}
We exhibit a connected graph on $24$ vertices with maximum degree $3$,
independence number $9$ and zero forcing number $11$, refuting a 2017
conjecture of \textsc{TxGraffiti} recorded as Conjecture~2 of the survey of
Davila, Brimkov and Pepper. The same construction with a different gadget
gives a connected \emph{cubic} graph on $36$ vertices with independence
number $15$ and zero forcing number $17$; the conjecture therefore fails also
in the cubic form in which the survey's Lean~4 appendix states it. In
particular $\Zf \le \alpha+1$ is not a universal bound for connected cubic
graphs, and the value $\Zf = \alpha+2$ is attained.
\end{abstract}

\maketitle

\section{The conjecture}

Colour a set $B \subseteq V(G)$ of vertices blue and the rest white. If a blue
vertex has exactly one white neighbour, it \emph{forces} that neighbour blue;
a vertex forces at most once. Iterating until no force is available yields the
\emph{closure} of $B$, which does not depend on the order of forces. The set
$B$ is a \emph{zero forcing set} if its closure is $V(G)$, and $\Zf(G)$ is the
least size of one. We write $\alpha(G)$ for the independence number.

\begin{conjecture}[\textsc{TxGraffiti}, 2017; {\cite[Conjecture~2]{Davila2025}}]
\label{conj:main}
If $G \not\cong K_4$ is a connected graph with $\Delta(G) \le 3$, then
$\Zf(G) \le \alpha(G)+1$.
\end{conjecture}

The text of \cite{Davila2025} states the conjecture for $\Delta(G) \le 3$, as
above, but the Lean~4 formalisation in its Appendix~A carries the additional
hypothesis \texttt{max\_degree G = min\_degree G} together with
\texttt{max\_degree G = 3}, that is, it asks that $G$ be cubic. The cubic graph
$G$ below refutes the conjecture under either reading; the subcubic graph $H$
refutes it as stated in the text, and is the smaller of the two.

Davila and Henning \cite{DavilaHenning2018} proved this for claw-free cubic
graphs, and Schuerger, Warnberg and Young \cite{SWY2024} proved that almost all
cubic graphs satisfy $\Zf \le \alpha+2$, exhibiting an infinite family with
$\Zf = \alpha+1$.

\section{The construction}

\begin{figure}[ht]
\centering
\begin{tikzpicture}[scale=1.05]
  \foreach \ang in {90,210,330}{
    \node[cvtx] (c\ang) at (\ang:0.75) {};
    \node[rvtx] (r\ang) at (\ang:1.85) {};
  }
  \draw[thick] (c90)--(c210)--(c330)--(c90);
  \foreach \ang in {90,210,330}{ \draw[thick] (c\ang)--(r\ang); }
  \foreach \ang in {90,210,330}{
    \begin{scope}[shift={(\ang:1.85)},rotate=\ang-90]
      \node[blob] (L\ang) at (-0.62,0.95) {\scriptsize $F$};
      \node[blob] (R\ang) at ( 0.62,0.95) {\scriptsize $F$};
      \draw (0,0)--(L\ang);
      \draw (0,0)--(R\ang);
    \end{scope}
  }
  \node[anchor=west] at (0.16,0.78) {$c_0$};
  \node[anchor=west] at (0.16,1.90) {$r_0$};
  \begin{scope}[shift={(4.7,0.35)}]
    \node[vtx] (z) at (0,0) {}; \node[below=1.5pt] at (z) {$z$};
    \node[vtx] (a) at (-0.42,0.78) {}; \node[left=1.5pt] at (a) {$a$};
    \node[vtx] (b) at (0.42,0.78) {}; \node[right=1.5pt] at (b) {$b$};
    \draw (z)--(a) (z)--(b) (a)--(b);
    \node at (0,-0.55) {$F = K_3$};
  \end{scope}
  \begin{scope}[shift={(6.9,0.35)}]
    \node[vtx] (z2) at (0,0) {}; \node[below=1.5pt] at (z2) {$z$};
    \node[vtx] (a2) at (-0.52,0.72) {}; \node[left=1.5pt] at (a2) {$a$};
    \node[vtx] (b2) at (0.52,0.72) {}; \node[right=1.5pt] at (b2) {$b$};
    \node[vtx] (p2) at (-0.36,1.52) {}; \node[left=1.5pt] at (p2) {$p$};
    \node[vtx] (q2) at (0.36,1.52) {}; \node[right=1.5pt] at (q2) {$q$};
    \draw (z2)--(a2) (z2)--(b2) (a2)--(p2) (a2)--(q2)
          (b2)--(p2) (b2)--(q2) (p2)--(q2);
    \node at (0,-0.55) {$F = D$};
  \end{scope}
\end{tikzpicture}
\caption{Left: the skeleton of $G(F)$ --- a triangle $c_0c_1c_2$ (black),
three vertices $r_0,r_1,r_2$ (grey), and six copies of the gadget $F$, each
joined to its $r_j$ at the attachment vertex. Right: the two gadgets used,
giving $H = G(K_3)$ on $24$ vertices and $G = G(D)$ on $36$ vertices.}
\label{fig:G}
\end{figure}
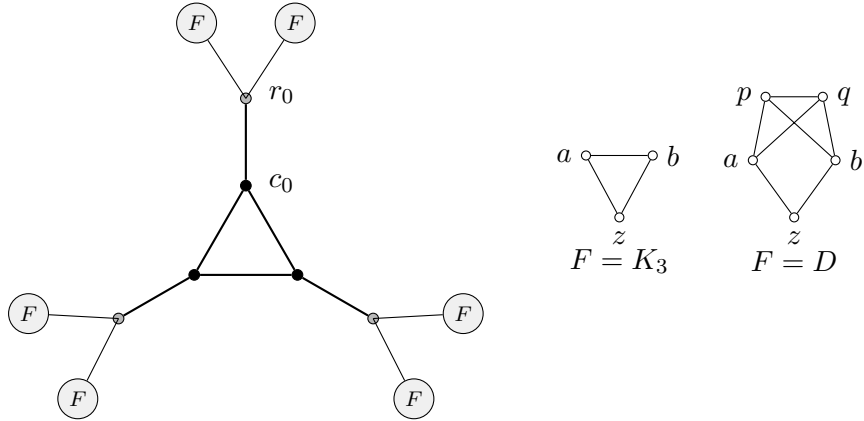

Call a graph $F$ with a distinguished vertex $z$ of degree $2$, all of whose
other vertices have degree $2$ or $3$, a \emph{gadget} with \emph{attachment
vertex} $z$. Given a gadget $F$, let $G(F)$ be the graph obtained from a
triangle $c_0c_1c_2$ by adding three vertices $r_0,r_1,r_2$ with edges
$c_jr_j$, and joining each $r_j$ to the attachment vertices of two disjoint
copies of $F$ (Figure~\ref{fig:G}). Then
\[
  |V(G(F))| = \underbrace{3}_{c_j} + \underbrace{3}_{r_j}
            + 6\,|V(F)| ,
\]
every $c_j$ and every $r_j$ has degree $3$, and each attachment vertex has
degree $3$ in $G(F)$. Hence $\Delta(G(F)) \le 3$, with equality throughout ---
that is, $G(F)$ is cubic --- exactly when every vertex of $F$ other than $z$
has degree $3$ in $F$. We call $r_j$ together with the two copies of $F$
attached to it the $j$th \emph{branch} $B_j$, so $|B_j| = 1 + 2|V(F)|$ and
$V(G(F))$ is the disjoint union of $B_0$, $B_1$, $B_2$ and the triangle. Each
edge $c_jr_j$ is a bridge, and each $r_j$ is the centre of a claw.

We use two gadgets.

\medskip
\noindent
\textbf{(i)} $F = K_3$, with vertices $z,a,b$ and $z$ the attachment vertex.
Write $H := G(K_3)$, a connected graph on $24$ vertices and $30$ edges with
$\Delta(H)=3$; the twelve vertices $a,b$ of the six copies have degree $2$, so
$H$ is subcubic but not cubic. Its \texttt{graph6} string is
\begin{center}\small
\verb|W{CGW_@?Y??@?@?@_@??@??K_????G??C??B??@????_??B|
\end{center}

\medskip
\noindent
\textbf{(ii)} $F = D$, the graph on $\{z,a,b,p,q\}$ with edges
$za, zb, ap, aq, bp, bq, pq$ --- that is, $K_{2,3}$ with parts $\{a,b\}$ and
$\{z,p,q\}$ plus the edge $pq$. Every vertex of $D$ other than $z$ has degree
$3$, so $D$ is the cubic completion of the triangle in the sense above.
Write $G := G(D)$, a connected \emph{cubic} graph on $36$ vertices and $54$
edges. Its \texttt{graph6} string is
\begin{center}\small
\verb|c{CGOKFC??_A?B?Bg???@??C??O??W??[??_???G??@???B???FG|\\
\verb|???????G???@????G????W????w???A?????@?????O????@_????F|
\end{center}
(the two lines concatenated). Both strings are also supplied as the files
\texttt{H24.g6} and \texttt{G36.g6}.

\section{The computation}

\begin{theorem}\label{thm:main}
$\alpha(H) = 9$ and $\Zf(H) = 11$; and $\alpha(G) = 15$ and $\Zf(G) = 17$.
Hence $\Zf = \alpha+2$ for both, and Conjecture~\ref{conj:main} is false ---
already on $24$ vertices, and also for connected cubic graphs. In particular
$\Zf \le \alpha+1$ is not a universal bound for connected cubic graphs.
\end{theorem}

All four quantities are computed exactly by finite search, and all four
computations are small enough to be repeated by hand-written code --- in
seconds in a compiled language, in a few minutes with the accompanying Python
script.

For the independence numbers, branch and bound over the vertex set terminates
immediately. A maximum independent set is obtained branch by branch: for $H$,
take $r_j$ together with one vertex other than $z$ from each of that branch's
two copies of $K_3$, three vertices per branch and $9$ in all; for $G$, take
$r_j$ together with the pair $a,b$ from each of its two copies of $D$, five
vertices per branch and $15$ in all.

For the upper bounds, an explicit zero forcing set is exhibited in each case
and its closure computed. For $H$, take $z$ and $a$ from five of the six
copies of $K_3$ and $a$ alone from the sixth, giving $5 \cdot 2 + 1 = 11$
vertices. For $G$, take $z$, $a$ and $p$ from five of the six copies of $D$
and $a,p$ alone from the sixth, giving $5 \cdot 3 + 2 = 17$ vertices.

For the lower bounds we use forts. A \emph{fort} is a nonempty
$F \subseteq V(G)$ such that no vertex outside $F$ has exactly one neighbour in
$F$; a set is a zero forcing set if and only if it meets every fort
\cite{Fast2018}. Let $t$ be the number of forts contained in a single branch
and $h$ the least size of a subset of a branch meeting all of them; by symmetry
these do not depend on $j$. Then every zero forcing set contains at least $h$
vertices of each branch, hence at least $3h$ in total, and the sets of size at
most $3h+1$ that survive this restriction number
\[
  m_h^{\,3} \;+\; 3\,m_h^{\,3} \;+\; 3\,m_{h+1}m_h^{\,2},
\]
where $m_h$ and $m_{h+1}$ count the subsets of a branch of size $h$ and $h+1$
meeting all its forts; the three terms count $|B| = 3h$, then $|B| = 3h+1$ with
the extra vertex on the triangle, and $|B| = 3h+1$ with the extra vertex in a
branch. The values are
\[
\begin{array}{c|cccccc}
 & |B_j| & t & h & m_h & m_{h+1} & \text{candidates} \\ \hline
 H & 7 & 12 & 3 & 12 & 25 & 17{,}712 \\
 G & 11 & 156 & 5 & 64 & 168 & 3{,}112{,}960
\end{array}
\]
Each of the $17{,}712$ sets of size at most $10$ in $H$, and each of the
$3{,}112{,}960$ sets of size at most $16$ in $G$, is checked directly and
fails to force. Hence $\Zf(H) \ge 11$ and $\Zf(G) \ge 17$.

A self-contained Python script using only the standard library,
\texttt{verify\_note.py}, accompanies this note; it rebuilds $H$ and $G$ from
the description above, verifies that each is connected with the stated degree
sequence and \texttt{graph6} string, and certifies all four of
$\alpha(H)=9$, $\Zf(H)=11$, $\alpha(G)=15$, $\Zf(G)=17$, printing one line per
claim.

The infinite family of which $G$ is the smallest cubic member, the minimum
order of a counterexample, and a repaired conjecture for bridgeless cubic
graphs are treated in a companion paper.

\end{document}